\documentclass[12pt,leqno]{article}
\usepackage{amsmath, amsthm, amsfonts, amssymb, color}
\usepackage{mathrsfs}
\newtheorem{thm}{Theorem}[section]

\theoremstyle{definition}

\newcommand{\scr}[1]{\mathscr #1}
\definecolor{wco}{rgb}{0.5,0.2,0.3}

\numberwithin{equation}{section} \theoremstyle{remark}

\newcommand{\ua}{\uparrow}

\title{{\bf  Exponential Convergence Rates of Second Quantization Semigroups and Applications }
\footnote{Supported in
 part by WIMCS and SRFDP.}
}
\author{
{\bf Chang-Song Deng$^{a)}$ and Feng-Yu Wang$^{a),b)}$\footnote{Corresponding author: wangfy@bnu.edu.cn,
F.Y.Wang@swansea.ac.uk \newline School of Mathematical Sciences, Beijing Normal University, Beijing 100875, China.}}\\
\footnotesize{$^{a)}$School of Mathematical Sciences, Beijing Normal
University, Beijing 100875, China}\\
 \footnotesize{$^{b)}$Department of Mathematics,
Swansea University, Singleton Park, SA2 8PP, UK} }
\begin{document}

\def\R{\mathbb R} \def\Z{\mathbb Z} \def\ff{\frac} \def\ss{\sqrt}
\def\N{\mathbb N} \def\kk{\kappa} \def\sgn{{\rm sgn}}\def\supp{{\rm
supp}}
\def\dd{\delta} \def\DD{\Delta} \def\vv{\varepsilon} \def\rr{\rho}
\def\<{\langle} \def\>{\rangle} \def\GG{\Gamma} \def\gg{\gamma}
\def\ll{\lambda} \def\LL{\Lambda} \def\nn{\nabla} \def\pp{\partial}
\def\d{\text{\rm{d}}} \def\bb{\beta} \def\aa{\alpha} \def\D{\scr D}
\def\E{\scr E} \def\si{\sigma} \def\ess{\text{\rm{ess}}}
\def\beg{\begin} \def\beq{\begin{equation}}  \def\F{\scr F}
\def\Ric{\text{\rm{Ric}}} \def\Hess{\text{\rm{Hess}}}\def\B{\scr B}
\def\e{\text{\rm{e}}} \def\ua{\underline a} \def\OO{\Omega} \def\b{\mathbf b}
\def\oo{\omega}     \def\tt{\tilde} \def\Ric{\text{\rm{Ric}}}
\def\cut{\text{\rm{cut}}} \def\P{\mathbb P} \def\ifn{I_n(f^{\bigotimes n})}
\def\fff{f(x_1)\dots f(x_n)} \def\ifm{I_m(g^{\bigotimes m})} \def\ee{\varepsilon}
\def\pm{\pi_{\mu}}   \def\p{\mathbf{p}}   \def\ml{\mathbf{L}}
 \def\C{\scr C}      \def\aaa{\mathbf{r}}
   \def\r{r} \def\BB{\scr B} \def\ttpsi{\tt\psi}
\def\gap{\text{\rm{gap}}} \def\prr{\pi_{\mu,\varrho}}  \def\r{\mathbf r}
\def\Z{\mathbb Z} \def\vrr{\varrho} \def\Sect{{\rm Sect}}\def\ent{{\rm Ent}}
\def\ii{{\rm i}_{\pp M}} \def\Ent{{\rm Ent}} \def\Var{{\rm Var}}
\def\q{\textbf{q}} \def\ess{{\rm ess}}\def\Eg{\scr E^\psi} \def\Etg{\tt{\scr E}^{\tt \psi}}

\maketitle
\begin{abstract}

\ \newline  Exponential
convergence rates in the $L^2$-tail norm and entropy are characterized for
the second quantization semigroups by using the corresponding base
Dirichlet form. This supplements the well known result on the
$L^2$-exponential convergence rate of second quantization
semigroups. As applications, birth-death type processes on Poisson
spaces and the path space of L\'evy processes are investigated.
\end{abstract}

\

\noindent
 AMS subject Classification:\ 60J75, 60H07, 60F10 \  \\
\noindent
 Keywords: Convergence rate, tail norm, entropy, second quantization,
  L\'evy processes.

 \vskip 2cm

\section{Introduction}

Let $E$ be a Polish space with Borel $\si$-field $\F$. Let $\mu$ be a non-trivial $\sigma$-finite measure on $(E,\F)$.
Let $(\E_0,\D(\E_0))$ be
 a symmetric Dirichlet form  on $L^2(\mu)$. Consider the
 configuration space
 \begin{eqnarray*}
\begin{aligned}
\GG:=\Big\{\gg=\sum_i \dd_{x_i}\ &\text{(at\ most\ countable)}:\
 x_i\in E\Big\},
\end{aligned}
\end{eqnarray*}
 where $\dd_x$ is the Dirac measure at $x$ and $\sum_\emptyset$ is
 regarded as the zero measure $0$ on $E$. Let $\F_\GG$ be the
 $\si$-field induced by $\{\gg\mapsto \gg(A): \ A\in \F\}$.
The
 Poisson measure with intensity $\mu$, denoted by $\pi_\mu$,  is the
 unique probability measure on $(\GG,\F_\GG)$ such that for any
 disjoint sets $A_1,\cdots, A_n\in \F$ with $\mu(A_i)<\infty$, $1\le i
\le n$,
$$\pi_\mu\big(\{\gg\in \GG: \ \gg(A_i)= k_i, 1\le i\le n\}\big)=
\prod_{i=1}^n \e^{-\mu(A_i)} \ff{\mu(A_i)^{k_i}}{k_i!},\ \ \
k_i\in \Z_+,\ 1\le i\le n.$$ This measure has the Laplace
transform
\beq\label{1.1} \pi_\mu(\e^{\<\cdot, f\>}) = \exp\left[\mu(\e^f-1)\right],\ \ \
f\in L^1(\mu)\cap L^\infty(\mu),\end{equation} where $\<\gg, f\>:= \gg(f)=\int_E f\,\d \gg.$

\medskip

\medskip

The second quantization of $(\E_0,\D(\E_0))$ is a symmetric
conservative Dirichlet form on $L^2(\pi_\mu)$ given by (see
e.g. \cite[Lemma 6.3]{Wu})
\beg{equation*}\beg{split} &\D(\E):= \Big\{F\in L^2(\pi_\mu):\ D_\cdot
F(\gg):= F(\gg+\dd_\cdot)- F(\gg)\in \D_e(\E_0),\
\pi_\mu\text{-a.e.}\ \gg, \\
&\qquad\qquad\qquad\qquad\qquad\qquad\qquad\qquad\qquad\qquad\E_0(D_\cdot
F, D_\cdot F)\in
L^1(\pi_\mu)\Big\},\\
&\E(F,G):= \int_{\GG}\E_0(D_\cdot F(\gg), D_\cdot
G(\gg))\pi_\mu(\d\gg),\ \ \ \ F,G\in
\D(\E),\end{split}\end{equation*} where $\D_e(\E_0)$ is the extended
domain of $\E_0$ (see \cite{f}).

\medskip

Let $P_t^0$ and $P_t$ be the semigroups associated to
$(\E_0,\D(\E_0))$ on $L^2(\mu)$ and $(\E,\D(\E))$ on $L^2(\pi_\mu)$
respectively. We aim to investigate the convergence rate of $P_t$ to
$\pi_\mu$ as $t\to\infty$ by using properties of the base Dirichlet
form.

\medskip

We would like to consider the following three kinds of exponential
convergence rates: \beg{enumerate} \item[(1)]  {\bf Exponential
convergence in the $L^2$-norm:} let $\ll_L$ be the largest   constant
such that
$$\| P_t - \pi_\mu \|_{L^2(\pi_\mu)\rightarrow L^2(\pi_\mu)}\le \e^{-\ll_L t},\ \ \ t\ge 0, $$
where $\pi_\mu$ is regarded as a linear operator from $L^2(\pi_\mu)$
to $\R$ by letting $\pi_\mu(F)=\int_\GG F\d\pi_\mu$.
\item[(2)]  {\bf Exponential convergence in the $L^2$-tail norm:} let $\ll_T$ be the largest   constant such that
$$\| P_t \|_{T}:=  \lim_{n\to\infty} \sup_{\pi_\mu(F^2)\le 1}
\|  1_{\{|P_t F|\ge n\}}P_t F\|_{L^2(\pi_\mu)}\le \e^{-\ll_T t},\ \ \
t\ge 0.$$ \item[(3)] {\bf Exponential convergence in entropy:} let
$\ll_E$ be the largest   constant such that
$$\pi_\mu ((P_t F)\log P_t F)\le \pi_\mu(F\log F)\e^{-\ll_E t}, \ \ \ t\ge 0, F\ge 0, \pi_\mu(F)=1.$$
\end{enumerate}

The exponential convergence rate in the $L^2$-norm is already well described
by the exponential decay rate of $P_t^0$, i.e.  (see \cite{Simon})
\beq\label{L}\ll_L= \ll_{L,0}:= \inf\big\{\E_0(f,f):\ f\in \D(\E_0),
\mu(f^2)=1\big\}.\end{equation} It is well known that $\ll_{L,0}$ is
the largest number such that
$$\|P_t^0f\|_{L^2(\mu)}\le \|f\|_{L^2(\mu)}\e^{-\ll_{L,0}t},\ \ \ t\ge
0, f\in L^2(\mu)$$ holds. See \cite{RW} and \cite{Wu} for a
criterion of the weak Poincar\'e inequality for  second quantization
Dirichlet forms.

\medskip

Due to the above fact,  in this paper  we will only consider $\ll_T$
and $\ll_E$. To study these two quantities, we first describe them
by using the Dirichlet form.

\medskip

\medskip

Since $\pi_\mu$ is a probability measure, by \cite[Theorem 3.3]{W02}
for $\phi\equiv 1$ we conclude that $\ll_T$ is the largest number
such that for any $C_1>\ll_T^{-1}$ the defective Poincar\'e
inequality
$$ \pi_\mu(F^2)\le C_1\E(F,F)+C_2\pi_\mu(|F|)^2,\ \ \
F\in\D(\E)$$ holds for some constant $C_2>0$. Consequently,
\beq\label{T} \ll_T=\lim_{n\to\infty}\inf\big\{\E(F,F)+n
\pi_\mu(|F|)^2:\ F\in\D(\E), \pi_\mu(F^2)=1\big\}.\end{equation} The
quantity $\ll_T$ is also related to the essential spectrum
$\si_{\ess}(\scr L)$ of the generator $\scr L$ associated to
$(\E,\D(\E))$. Precisely, we have $$\ll_T\ge \inf\si_{\ess}(-\scr
L)$$ and the equality holds provided for some $t>0$ the operator
$P_t$ has an asymptotic density w.r.t. $\pi_\mu$ (see \cite[Theorem
3.2.2]{w05}).

\medskip

\medskip

Next, it is easy to check that  $\ll_E$ is the largest number such
that the $L^1$ log-Sobolev inequality
\begin{eqnarray}
\begin{aligned}\label{entropy}
{\rm Ent}_{\pi_\mu}(F)&:= \pi_\mu(F\log F)-\pi_\mu(F)\log\pi_\mu(F)\\
&\le \ff 1 {\ll_E} \E(F,\log F),\ \ \
F\in\D(\E),\inf F>0
\end{aligned}
\end{eqnarray}
holds. That is (see \cite[Theorem 1.1]{MZ}),
\beq\label{E} \ll_E= \inf\left\{\ff{\E(F,\log F)}{{\rm
Ent}_{\pi_\mu}(F)}:\ \inf F>0, F\in \D(\E),
{\Ent}_{\pi_\mu}(F)>0\right\}.\end{equation}
We remark that for $F\in\D(\E)$ with $\inf F>0$, one has $\log F\in\D(\E)$ so that
$\E(F,\log F)$ exists.

\medskip

\medskip

Finally, we would like to mention that  the   log-Sobolev inequality
introduced in \cite{G}
\beq\label{LS} \Ent_{\pi_\mu} (F^2) \le C \E(F,F),\ \ \ F\in
\D(\E) \end{equation} for some constant $C>0$ implies that
$\ll_E\ge 4/C$ (see e.g. \cite[Theorem 1.2]{MZ}). But it is easy to see that
the second quantization Dirichlet form does not satisfy the
log-Sobolev inequality (see \cite{S} and the first page of
\cite{wu00}). Indeed, given nonnegative function $f\in
L^\infty(\mu)\cap L^1(\mu)\cap \D(\E_0)$, applying (\ref{LS}) to
$F(\gg):=\e^{\gg(f)}$ and using (\ref{1.1}) we obtain
$$\int_E(2f\e^{2f}-\e^{2f}+1)\,\d\mu \le C\E_0\big(\e^{f}-1,\e^f-1\big).$$
Replacing $f$ by $\log (n f+1)$ which is once again in
$L^\infty(\mu)\cap L^1(\mu)\cap \D(\E_0)$, we obtain
$$\ff1{n^2\log n}\int_E\big\{2(nf+1)^2\log (nf+1) -
(nf+1)^2+1\big\}\,\d\mu \le \ff C {\log n} \E_0(f,f).$$ Letting
$n\to\infty$ we arrive at $\mu(f^2)\le 0$ which is impossible if
$f$ is non-trivial.

\

It is now the place to  state our main result of the paper where
$\ll_E$ and $\ll_T$ are described by using the base Dirichlet form
$(\E_0, \D(\E_0)).$

\beg{thm}\label{T1.1} We have
\beq\label{LE} \ll_E=
\inf\left\{\ff{\E_0\big(\e^f-1,f\big)}{\mu(f\e^f-\e^f+1)}:\ f\in\D(\E_0)\cap
L^\infty(\mu), \mu(f^2)>0\right\} \end{equation} and
\beq\label{LT} \ll_{L,0}\le \ll_{T}\le \ll_{T,0}:= \lim_{n\to\infty}
\inf\big\{\E_0(f,f)+n\mu(|f|)^2:\ f\in\D(\E_0),
\mu(f^2)=1\big\}.\end{equation} \end{thm}

\

To derive the exact value of these two quantities, let us decompose
the Dirichlet form $\E_0$ into three parts: the diffusion part, the
jump part and the killing part. We will see in the next result that
in many cases $\ll_E$ is determined merely by the killing term.

\

Let $W$ be a nonnegative measurable function on $E$, $\scr A\subset L^1(W\mu)\cap L^\infty(\mu)$ be a linear
subspace,   $q \ge 0$ be a symmetric measurable function on $E\times
E$, and $\GG_1: \scr A\times \scr A\to L^1(\mu)$ be a nonnegative
definite  bilinear map    such that
 \beg{enumerate}
\item[(i)]  $\scr A$ is dense in $L^2\big((1+W)\mu\big)$;\item[(ii)]  If $f\in
\scr A$ and $\phi: \R\to \R$ is Lipschitz continuous  with
$\phi(0)=0$, then $\phi(f)\in \scr A$;\item[(iii)]  For any $f\in
\scr A,\ \int_{E\times E} |f(x)-f(y)|^2 q(x,y)\mu(\d x)\mu(\d
y)<\infty;$ \item[(iv)]   $\GG_1(f,\phi(g))= \phi'(g) \GG_1(f,g)$ holds
for any $\phi\in C^1(\R)$ with $\phi(0)=0$ and any $f,g\in\scr A$.
\end{enumerate}

\medskip

Consider the following diffusion-jump type quadric form with
potential:
\beq\label{2331} \beg{split}\E_0(f,g)&:= \mu\big(\GG_1(f,g)+ Wfg \big) \\
&\quad+ \ff 1 2 \int_{E\times E} (f(x)-f(y))(g(x)-g(y)) q(x,y)\mu(\d
x)\mu(\d y),\ \ \ f,g\in \scr A.\end{split}\end{equation} Assume that
$(\E_0,\scr A)$ is closable such that its closure $(\E_0,\D(\E_0))$ is
a Dirichlet form on $L^2(\mu)$. When $\GG_1=0$, $q=0$ and
$W\equiv 1$, the framework goes back to \cite{wu00} where the
Poincar\'e inequality and the $L^1$ log-Sobolev inequality with
constant 1 are proved. The contribution of our next result is to
confirm that these inequalities    are sharp under a more general
framework.

\beg{cor}\label{C1.2}  Let $(\E_0,\D(\E_0))$ be given in
$(\ref{2331})$ such that {\rm (i)--(iv)} hold.
\beg{enumerate}\item[$(1)$]\ If there exists a sequence of
nonnegative functions $\{f_n\}_{n\geq1}\subset \scr A$ such that
$\{f_n>0\}\uparrow E$ as $n\uparrow\infty$, then $\ll_E=\ess_\mu\inf
W.$
\item[$(2)$]\  Let $\GG_1=0$ and $q=0$, and let $\mu$ be finite on bounded sets. If $\supp\mu\cap \{W<\vv\}$ is
uncountable whenever $\mu(W<\vv)>0\ ($it is the case if $\mu$ does
not have atom$)$, then $\ll_L=\ll_T=\ess_\mu\inf W.$
\end{enumerate}\end{cor}

To conclude this section, we present below an example to illustrate
Corollary \ref{C1.2}(1).

\paragraph{Example 1.1.} Let $E$ be a connected (not necessarily
complete) Riemannian manifold and $V$ a locally bounded measurable
function. Let $\mu(\d x)=\e^{V(x)}\d x$ with $\d x$ the volume
measure. Then we take $\scr A$ to be the set of all Lipschitz
continuous functions on $E$ with compact supports. It is trivial
that conditions (i) and (ii) hold and $\scr A\subset
L^1(W\mu)\cap L^\infty(\mu)$ provided $W$ is locally bounded.
Define
$$\GG_1(f,g)=\<\nn f,\nn g\>,\ \ \ f,g\in \scr A.$$
Then condition (iv) holds. Finally, let $\rr(x,y)$ be the Riemannian
distance between $x$ and $y$. If $q(x,y)$ satisfies
\beq\label{2334}\int_{K\times E} \big(\rr(x,y)^2\land 1\big)q(x,y)\mu(\d
x)\mu(\d y)<\infty\end{equation} for any compact subset $K$ of $E$,
then (iii) is satisfied. Thus, by Corollary \ref{C1.2}(1) where the
required sequence $\{f_n\}_{n\geq1}$ automatically exists according to the
definition of $\scr A$, we have
$$\ll_E=\ess_\mu\inf W.$$  In particular, let $\mu$ be the Lebesgue measure and   $E$   a   bounded open
domain   in $\R^d$ (it is complete under a compatible  metric), a typical choice of $q(x,y)$ such that (\ref{2334}) holds is $\ff
1 {|x-y|^{\aa+d-1}}$ for $\aa\in [0,2)$. Moreover, if $E=\R^d$ and $\mu(\d x)=\d x$, then (\ref{2334}) holds
for this   $q(x,y)$ with $\aa\in (1,2)$.

\

The remainder of the paper is organized as follows. In Section 2   complete proofs of Theorem \ref{T1.1} and Corollary
\ref{C1.2} are presented;    In Section 3 the exponential convergence rates  are
considered for birth-death type Dirichlet forms on
$L^2(\pi_\mu)$ with a weighted function on $\GG\times E$;  and  in Section 4
results derived in Section 3 are applied to    the path space
of L\'evy processes by following the line of \cite{wu00}.

\section{Proofs of Theorem \ref{T1.1} and Corollary \ref{C1.2}}

\beg{proof}[Proof of (\ref{LE})] We first remark that for any $f\in
\D(\E_0)\cap L^\infty(\mu)$ one has $\e^f-1\in \D(\E_0)$, since the
function $\phi(r):= \e^r-1$ is locally Lipschitz continuous and
$\phi(0)=0.$  Therefore, it suffices to show that for any $\ll>0$,
the $L^1$ log-Sobolev inequality
\beq\label{LS1} {\rm Ent}_{\pi_\mu}(F) \le \ff 1\ll  \E(F,\log F),\ \ \
F\in\D(\E),\inf F>0\end{equation} is equivalent to  \beq\label{LS2}
\mu\big(f\e^f -\e^f+1\big)\le \ff 1 \ll \E_0\big(\e^f-1, f\big),\ \ \
f\in \D(\E_0)\cap L^\infty(\mu).\end{equation}

\medskip

(a) (\ref{LS2}) implies (\ref{LS1}).  It suffices to prove
(\ref{LS1}) for $F\in \D(\E)\cap L^\infty(\pi_\mu)$ with $\inf F>0$.
In this case we have $g_\gg:= \ff{F(\gg+\dd_\cdot)}{F(\gg)} -1\in
\D_e(\E_0)$  for $\pi_\mu$-a.e. $\gg\in\GG$. Since $\ff{\sup F}{\inf F}\ge g_\gg+1>0$, it follows that
$$f_\gg:= \log (g_\gg +1)\in \D_e(\E_0)\cap L^\infty(\mu)$$
for $\pi_\mu$-a.e. $\gg\in\GG$. By (\ref{LS2}) which holds also for $f\in
\D_e(\E_0)\cap L^\infty(\mu),$ we have
\beq\label{LS2'}\ll\int_E(f_\gg\e^{f_\gg}- \e^{f_\gg}+1)\,\d\mu \le
\E_0\big(\e^{f_\gg}-1, f_\gg\big) =\E_0\big( g_\gg, \log(g_\gg
+1)\big).\end{equation} On the other hand, by  the modified log-Sobolev
inequality presented in \cite[Theorem 1.1]{wu00} (note that
$\Phi(r)=r\log r$ therein), it holds that
\beq\label{MLS} {\rm Ent}_{\pi_\mu} (F)\le \int_{\GG}\pi_\mu(\d\gg)
\int_E \big\{D_z(F\log F)(\gg)   -\big(1+\log F(\gg)\big)D_z
F(\gg)\big\} \mu(\d z).\end{equation} Since
$$D_z F(\gg)= F(\gg) \big(\e^{f_\gg}(z) -1\big),\ \ \ \log\ff {F(\gg+\dd_z)}{F(\gg)} =f_\gg(z),$$
it is not hard to verify that
\beg{equation*}\beg{split} &D_z(F\log F)(\gg)   -\big(1+\log
F(\gg)\big)D_z F(\gg)= F(\gg+\dd_z) \log
\ff{F(\gg+\dd_z)}{F(\gg)}  -D_z F(\gg)\\
&= (D_z F(\gg))\Big(\log \ff{F(\gg+\dd_z)}{F(\gg)} -1\Big) + F(\gg)
\log \ff{F(\gg+\dd_z)}{F(\gg)}\\
&= F(\gg)\big\{ \big(\e^{f_\gg}-1\big)(f_\gg-1) +f_\gg\big\}(z)
=F(\gg) \big(f_\gg
\e^{f_\gg}-\e^{f_\gg}+1\big)(z).\end{split}\end{equation*} Combining
this with (\ref{LS2'}) and (\ref{MLS}), we obtain
\beg{equation*}\beg{split}\ll\,\Ent_{\pi_\mu}(F) &\le \ll \int_\GG
F(\gg) \pi_\mu(\d\gg) \int_E (f_\gg \e^{f_\gg} -\e^{f_\gg}
+1)\d\mu\\
&\le \int_\GG F(\gg) \E_0(g_\gg, \log (g_\gg+1))\pi_\mu(\d\gg)\\
&= \int_\GG \E_0(D_\cdot F, D_\cdot \log F)\d\pi_\mu = \E(F,\log
F).\end{split}\end{equation*}

\medskip

(b) (\ref{LS1}) implies (\ref{LS2}). We first consider
$f\in \D(\E_0)\cap L^\infty(\mu)\cap L^1(\mu)$. Let $F(\gg)= \e^{\gg(f)}$. By (\ref{1.1})
we have $F\in L^2(\pi_\mu)$ and
\beq\label{2.2} \Ent_{\pi_\mu}(F)= \pi_\mu(F) \int_E(f\e^f-\e^f
+1)\,\d\mu.\end{equation} Moreover, for any $\vv>0$ one has $F+\vv \in
\D(\E),\ \inf (F+\vv)>0$ and
\begin{eqnarray}
\begin{aligned}\label{2.3} &\E\big(F+\vv, \log (F+\vv)\big) \\&\qquad=\int_\GG
F(\gg)\left\{\E_0\big(\e^f-1, f\big)+\E_0\Big(\e^f-1, \log
\ff{\e^{\gg(f)}+\vv
\e^{-f}}{\e^{\gg(f)}+\vv}\Big)\right\}\pi_\mu(\d\gg).
\end{aligned}
\end{eqnarray}
Since $\phi(s):= \log \ff{\e^{\gg(f)}+\vv \e^{-s}}{\e^{\gg(f)}+\vv}$
satisfies $\phi(0)=0$ and $|\phi'(s)|\le 1$, we get
\begin{eqnarray*}
\begin{aligned}
\left|\E_0\Big(\e^f-1, \log \ff{\e^{\gg(f)}+\vv
\e^{-f}}{\e^{\gg(f)}+\vv}\Big)\right| &\le \sqrt{\E_0\big(\e^f-1,\e^f-1\big)\E_0\big(\phi(f),\phi(f)\big)}\\
&\le\ss{\E_0\big(\e^f-1,\e^f-1\big)\E_0(f,f)}<\infty.
\end{aligned}
\end{eqnarray*}
Thus, by (\ref{2.3})
and the dominated convergence theorem we arrive at
\beq\label{2.4} \beg{split} &\lim_{\vv\downarrow 0} \E\big(F+\vv,
\log (F+\vv)\big) \\&=\int_\GG F(\gg) \E_0\big(\e^f-1,f\big)\pi_\mu(\d\gg)
+\int_\GG F(\gg) \lim_{\vv\downarrow 0} \E_0\Big(\e^f-1,
\log
\ff{\e^{\gg(f)}+\vv \e^{-f}}{\e^{\gg(f)}+\vv}\Big)\pi_\mu(\d\gg)\\
& =\pi_\mu(F) \E_0\big(\e^f-1, f\big).\end{split}\end{equation} Therefore,
first applying (\ref{LS1}) to $F+\vv$ then letting $\vv\downarrow0$, we
obtain (\ref{LS2}) from (\ref{2.2}) and (\ref{2.4}).

In general, for any  $f\in\D(\E_0)\cap L^\infty(\mu)$, let
$$f_n=\left(f-\frac1n\right)^+-\left(f+\frac1n\right)^-,\quad n\geq1.$$
Then it is easy to see that $f_n\in \D(\E_0)\cap L^\infty(\mu)\cap L^1(\mu)$ and $f_n\rightarrow f$ in
$\D(\E_0)\cap L^\infty(\mu)$. Therefore, (\ref{LS2}) holds.
\end{proof}

\medskip

\beg{proof}[Proof of (\ref{LT})] Since it is well known that
$$\ll_L= \inf\{\E(F,F): F\in \D(\E), \pi_\mu(F^2)-\pi_\mu(F)^2=1\},$$
   (\ref{L}) and (\ref{T})  imply $\ll_T\ge
\ll_{L,0}$. So, it remains to prove  $\ll_T\le \ll_{T,0}$. If
$0<\ll <\ll_T$, then there exists $C>0$ such that
\beq\label{DD*}\pi_\mu(F^2)\le \ff 1 \ll \E(F,F) + C\pi_\mu(F)^2,\ \ \ F\in \D(\E),
F\ge 0.\end{equation}  For any $f\in \D(\E_0)$, letting $F(\gg)= \gg(|f|)$
we have $\E(F,F)=\E_0(|f|,|f|)\le \E_0(f,f)$ and (see e.g.
\cite[Proof of Lemma 7.2]{RW})
$$ \pi_\mu(F^2)= \mu(f^2)+\mu(|f|)^2,\ \ \pi_\mu(F)= \mu(|f|).$$
Therefore, it follwos from (\ref{DD*}) that
$$\mu(f^2)\le \ff 1 \ll \E_0(f,f) + (C-1)\mu(|f|)^2,\ \ \ \  f\in \D(\E_0).$$
This implies that $\ll_{T,0}\ge \ll$ holds for any $\ll<\ll_T.$
Hence,  $\ll_T\le \ll_{T,0}$.
\end{proof}

\medskip

To prove Corollary \ref{C1.2}, we need the following fundamental lemma. We
include a simple proof for completeness.

\beg{lem}\label{L2.1} Let $\nu$ be a measure on $E$ such that $\nu$ is finite
on bounded sets. If there exists a constant $c>0$ such that
$\nu(f^2)\le c\nu(|f|)^2$ holds for all $f\in L^2(\nu)$, then $\supp\nu$ is at most countable. If moreover $\nu(E)<\infty$ then
$\supp\nu$ is finite.\end{lem}

\beg{proof} Since $\nu$ is finite on bounded sets and   $E$ is
separable, there exists a sequence of open sets $\{G_n\}_{n\ge 1}$
such that $\cup_{n\ge 1}G_n=E$ and $\nu(G_n)<\infty$ for $n\ge 1$. Now we fix $n\geq1$.
Suppose there are $m$ many different points $\{x_i\}_{i= 1}^m$
   in $\text{supp}\nu\cap G_n$, where $m\geq1$. For each $i$ there exists $r_i>0$ such that    $B_i:= \{x:
 d(x,x_i)<r_i\}\subset G_n$ and $\{B_i\}_{i=1}^m $ are disjoint.
 Since $x_i$ is in the support of $\nu$, we have $\nu(B_i)>0$ for
 each $i\in\{1,\cdots,m\}$. Moreover, since $$\sum_{i=1}^m \nu(B_i)=\nu\left(\bigcup_{i=1}^m B_i\right)\le
 \nu(G_n)<\infty,$$
there exists $i_0\in\{1,\cdots,m\}$ such that     $$0<\nu(B_{i_0})\le \ff 1 m \nu(G_n).$$ But
 applying $\nu(f^2)\le c\nu(|f|)^2$ to $f=1_{B_{i_0}}$ we obtain $\nu(B_{i_0})\ge
1/c$.  Therefore,  $m\le  c\nu(G_n)$. This means
 that for each fixed $n\ge 1$ the set $\supp\nu\cap G_n$ is   finite,  so that $\supp\nu$ is at most
 countable. The second assertion follows from the same argument by
 taking $G_n=E$.\end{proof}

\medskip

\beg{proof}[Proof of Corollary \ref{C1.2} (1)] Since for any $r\in
\R$ one has
$$r(\e^r-1)\ge r\e^r -\e^r +1,$$ it holds that
\beg{equation*}\beg{split}\E_0\big(\e^f -1, f\big) &\ge
\int_EWf(\e^f-1)\,\d\mu
\ge (\ess_\mu\inf W)\int_E f(\e^f-1)\,\d\mu\\
&\ge (\ess_\mu\inf W)
\int_E (f\e^f -\e^f+1)\,\d\mu.\end{split}\end{equation*} Therefore, it
follows from (\ref{LE}) that $\ll_E\ge \ess_\mu\inf W.$

\medskip

On the other hand, let $g\in\scr A$ be a fixed nonnegative function.
For any $n\ge 1$, applying (\ref{LE}) to $f:= 2\log (ng +1)\in \scr
A\subset \D(\E_0)\cap L^\infty(\mu)$ and noting that by (iv)
$$\GG_1\big((ng+1)^2-1,2\log(ng+1)\big)=4n^2\GG_1(g,g),$$
we obtain
\beg{equation*}\beg{split} &\ll_E \int_E\big\{(ng+1)^2\log\left[(ng+1)^2\right]
-(ng +1)^2 +1\big\}\,\d\mu\\ &\le \E_0\big((ng+1)^2-1, 2\log (ng+1)\big)\\
&=
\int_E\big\{4n^2\GG_1(g,g)+W(n^2g^2+2ng)\log\left[(ng+1)^2\right]\big\}\,\d\mu\\
&\quad + \int_{E\times E}
\big\{(ng(x)+1)^2-(ng(y)+1)^2\big\}\Big(\log \ff
{ng(x)+1}{ng(y)+1}\Big) q(x,y) \mu(\d x)\mu(\d
y).\end{split}\end{equation*} Multiplying both sides by $\ff 1
{n^2\log n}$ and letting $n\to\infty$, by the dominated convergence
theorem we arrive at
\beq\label{2332} 2\mu(g^2(\ll_E-W))\le
\limsup_{n\to\infty}\int_{E\times E} G_n(x,y) q(x,y)\mu(\d x)\mu(\d
y),\end{equation} where
\beg{equation*}\beg{split} 0\le G_n(x,y)&:=
\ff{(ng(x)+1)^2-(ng(y)+1)^2}{n^2\log n} \log
\ff{ng(x)+1}{ng(y)+1}\\
&\le \ff{(ng(x)+ng(y)+2)\log(n [g(x)\lor g(y)]+1)}{n\log n}
|g(x)-g(y)|\\
&\le c|g(x)-g(y)|\end{split}\end{equation*}
for $\mu$-a.e. $x,y\in E$ and some constant $c>0$
since $g\in L^\infty(\mu)$. Thus, by (iii) and the dominated convergence
theorem it follows that
\beg{equation*}\beg{split}&\lim_{n\to\infty}\int_{E\times E}
G_n(x,y) q(x,y)\mu(\d x)\mu(\d y)\\&=\int_{E\times E}\lim_{n\to\infty} G_n(x,y) q(x,y)\mu(\d x)\mu(\d y)\\
&=\int_{E\times E} \big(g(x)^2-g(y)^2\big) \big(1_{\{g>0\}}(x)- 1_{\{g>0\}}(y)\big)
q(x,y)\mu(\d x)\mu(\d y).\end{split}\end{equation*} Combining this
with (\ref{2332}) and using the symmetry of $q(x,y)$ we get
\beq\label{2333} \beg{split}&\mu(g^2(\ll_E-W))\\ &\le \ff12\int_{E\times E} \big(g(x)^2-g(y)^2\big) \big(1_{\{g>0\}}(x)-
1_{\{g>0\}}(y)\big) q(x,y)\mu(\d x)\mu(\d y)\\
&=\int_{\{g>0\}\times\{g=0\}}g(x)^2q(x,y)\mu(\d x)\mu(\d y),\ \ \ g\in\scr A, g\ge0.
\end{split}\end{equation}

Next, let $E_n=\{f_n>0\}$. For any $n,m\ge 1$, applying (\ref{2333})
to $g_{nm}:=g+f_n/m$ we have
\beg{equation*}\beg{split} &\mu(g^2_{nm}(\ll_E-W))\\&\le
\int_{(E_n\cup \{g>0\})\times (E_n^c\cap \{g=0\})}g_{nm}(x)^2
q(x,y)\mu(\d x)\mu(\d y)\\
&\le \Big(\|g\|_\infty+\ff{\|f_n\|_\infty}m\Big)\int_{
\{g>0\}\times
(E_n^c\cap \{g=0\})}\Bigg\{|g(x)-g(y)| \\
&\qquad\qquad\qquad\qquad\qquad\qquad\qquad\qquad+\ff 1 m
|f_n(x)-f_n(y)|\Bigg\}q(x,y)\mu(\d x)\mu(\d y)\\
&\quad + \ff 1 {m^2}\|f_n\|_\infty \int_{(E_n\setminus\{g>0\})\times
(E_n^c\cap \{g=0\})}|f_n(x)-f_n(y)| q(x,y)\mu(\d x)\mu(\d
y).\end{split}\end{equation*}It follows by letting $m\to\infty$ that
$$\mu(g^2(\ll_E-W))\le \|g\|_\infty \int_{\{g>0\}\times E_n^c}
|g(x)-g(y)| q(x,y)\mu(\d x)\mu(\d y).$$ Finally, letting
$n\to\infty$ we conclude that $\mu(g^2(\ll_E-W))\le 0$ for any
nonnegative $g\in \scr A$. Since $\phi(x)=|x|$ is Lipschitz continuous with
$\phi(0)=0$, it holds that $\mu(g^2(\ll_E-W))\le 0$ for any $g\in \scr A$.
Noting that $\scr A$ is dense in $L^2\big((1+W)\mu\big)$, then it is trivial to see that
$\ll_E\le\ess_\mu\inf W$. This completes the proof.
\end{proof}

\medskip

\medskip

\beg{proof}[Proof of Corollary \ref{C1.2} (2)] Let $\GG_1=0$ and $q=0$. Then $\E_0(f,g)=\mu(Wfg)$.
In this case, we have
$$\ll_{L,0}=\inf_{f\in  L^2(\mu),\mu(f^2)>0}\frac{\mu(Wf^2)}{\mu(f^2)}=\ess_\mu\inf W.$$
So, by Theorem \ref{T1.1}, it suffices
to show that $\ll_{T,0}\le \ess_\mu\inf W$. If
$\ll_{T,0}>\ess_\mu\inf W$ then there exist $0<r<\{\ess_\mu\inf
W\}^{-1}$ and $c>0$ such that
\beq\label{CC}\mu(f^2)\le r\E_0(f,f)+c\mu(|f|)^2= r \mu(Wf^2)
+c\mu(|f|)^2,\ \ \  f\in L^2(\mu)\end{equation} holds. Take $\vv\in
(0,r^{-1})$ such that $\mu(W<\vv)>0$. Let $\mu_\vv=
1_{\{W<\vv\}}\mu$. Using $f1_{\{W<\vv\}}$ to replace $f$, we obtain
from (\ref{CC}) that
$$\mu_\vv (f^2)\le \ff c {1-r\vv} \mu_\vv(|f|)^2,\ \ \ f\in L^2(\mu_\vv).$$ Thus, according to Lemma
\ref{L2.1} supp$\mu_\vv$ is at most countable. This is contradictive
to the assumption that supp$\mu\cap\{W<\vv\}$ is
uncountable.\end{proof}

\section{Birth-death type Dirichlet forms on $L^2(\pi_\mu)$}

Let $\psi$ be a nonnegative measurable function on $\GG\times E$ such that
$$\psi_\mu(z):=\int_\GG \psi(\gg,z)\pi_\mu(\d \gg)<\infty,\ \ \ \mu\text{-a.e.}\ z\in E.$$ Consider the quadric form
\beg{equation*}\beg{split} &\Eg(F,G) :=\int_{\GG\times E} \big(F(\gg+\dd_z)-F(\gg)\big)\big(G(\gg+\dd_z)-G(\gg)\big)
\psi(\gg,z)\pi_\mu(\d\gg)\mu(\d z),\\
&\D(\Eg):= \{F\in L^2(\pi_\mu):\
\Eg(F,F)<\infty\}.\end{split}\end{equation*} According to
Propositions \ref{P3.3} and \ref{P3.4} below, $(\Eg,\D(\Eg))$ is a
conservative symmetric Dirichlet form on $L^2(\pi_\mu)$, which is
regular provided $\mu(\psi_\mu)<\infty$. Obviously, if $\psi(\gg,z)$
does not depend on $\gg$, then $\E^\psi$ goes back to the second
quantization Dirichlet form for $\E_0(f,g):=\mu(\psi fg)$ with
$\D(\E_0)= L^2((1+\psi)\mu).$

\

\beg{thm}\label{T3.1} Let $\ll_{L}(\psi), \ll_{T}(\psi)$ and
$\ll_{E}(\psi)$ be, respectively, the exponential convergence rates
in the $L^2$-norm, the $L^2$-tail norm and entropy for the semigroup associated
to $(\Eg,\D(\Eg))$. \beg{enumerate} \item[$(1)$]
In general, we have $\ess_{\pi_\mu\times\mu}\inf \psi\le \ll_{L}(\psi), \ll_{E}(\psi)\le \ess_\mu\inf \psi_\mu$.
If $\psi(\gg,z)$ is independent of $\gg$, then
$\ll_{L}(\psi)=\ll_{E}(\psi)=\ess_\mu\inf \psi$.
\item[$(2)$] Let $\mu$ do not have atom and be finite on bounded sets. Then
 $\ess_{\pi_\mu\times\mu}\inf \psi\le \ll_{T}(\psi)  \le \ess_\mu\inf
\psi_\mu.$  If  moreover $\psi(\gg,z)$ does not depend on $\gg$, then
$\ll_{T}(\psi)=\ess_\mu\inf \psi$.   \end{enumerate}\end{thm}

\beg{proof} (1) Let $\E$ be the second quantization Dirichlet form for
$\E_0(f,g):= \big(\ess_{\pi_\mu\times\mu}\inf\psi\big)\mu(fg)$. Obviously, we
have $\Eg\ge\E$. Combining this with Corollary \ref{C1.2} and
(\ref{L}) we conclude that
$$\lambda_L(\psi)\wedge\lambda_E(\psi)\ge
\ess_{\pi_\mu\times\mu}\inf\psi.$$  Consequently, it suffices to prove the
desired upper bound estimate.

\medskip

Taking $F(\gg)= \gg(f)$ for nonnegative $f\in L^1(\mu)\cap
L^\infty(\mu)$, we see that the defective Poincar\'e inequality
\beq\label{DD} \pi_\mu(F^2)\le
C_1\Eg(F,F)+C_2\pi_\mu(F)^2\end{equation}
implies that
\beq\label{D2} \mu(f^2)\le C_1 \mu(\psi_\mu f^2)
+(C_2-1)\mu(f)^2.\end{equation} Thus, (\ref{DD}) for $C_2=1$ (i.e.
the Poincar\'e inequality) implies that $C_1\ge (\ess_\mu
\inf\psi_\mu)^{-1}.$ That is, $\ll_{L}(\psi)\le \ess_\mu\inf\psi_\mu.$

\medskip

On the other hand, according to (b) in the proof of (\ref{LE}), the $L^1$
log-Sobolev inequality
\beq\label{LL1} \pi_\mu(F\log F)\le \ll \Eg(F,\log F)+
\pi_\mu(F)\log\pi_\mu(F)\end{equation} for $F(\gg):=\e^{\gg(f)}$
implies that
$$\mu(f\e^f-\e^f+1)\le \ll \mu\big(\psi_\mu (\e^f-1)f\big),\ \ \ f\in
L^\infty(\mu)\cap L^1(\mu).$$  Hence, by the proof of Corollary
\ref{C1.2} for $ W=\psi_\mu$, $\GG_1=0$ and $q=0$, we conclude that
(\ref{LL1}) implies $\ll\ge (\ess_\mu\inf\psi_\mu)^{-1}$. This means
that $\ll_{E}(\psi)\le \ess_\mu\inf\psi_\mu.$

\medskip

(2) Assume that $\mu$ does not have atom and is finite on bounded sets. According to Theorem \ref{T1.1}, we obtain
$$\lambda_T\geq\lambda_{L,0}=\ess_{\pi_\mu\times\mu}\inf\psi.$$
Finally, by Lemma \ref{L2.1}, (\ref{D2})
for any $C_2>0$ implies that $C_1\ge (\ess_\mu\inf \psi_\mu)^{-1}.$
Now we conclude that $\ll_{T}(\psi)\le \ess_\mu\inf\psi_\mu$ and the proof is completed.
\end{proof}

\

The remainder of this section   devotes to characterizing the form
$(\Eg, \D(\Eg))$.  To see that it is a   Dirichlet form on $L^2(\pi_\mu)$, we
need the following quasi-invariant property of the map $\gg\mapsto
\gg+\dd_z$.

\beg{lem}\label{L3.2} If $A\in \F_\GG$ is a $\pi_\mu$-null set, then

$$\tt A:= \{(\gg,z)\in \GG\times E:\ \gg+\dd_z\in A\}$$
is a $(\pi_\mu\times\mu)$-null set. \end{lem}

\beg{proof} We shall make use of the Mecke identity \cite{me} (see
also \cite{ro}), i.e. \beq\label{M} \int_{\GG\times E} H(\gg+\dd_z,
z)\pi_\mu(\d\gg)\mu(\d z)= \int_{\GG\times E} H(\gg,z)\gg(\d z)\pi_\mu(\d\gg)
\end{equation} holds for any measurable function $H$ on
$\GG\times E$ such that one of the above integrals exists.
Applying (\ref{M}) to $H(\gg,z)= 1_A(\gg)$ and noting that $\pi_\mu(A)=0$, we obtain
\begin{eqnarray*}
\begin{aligned}
(\pi_\mu\times\mu)(\tt A)&=\int_{\GG\times E}1_A(\gg+\delta_z)\pi_\mu(\d\gg)\mu(\d z)\\
&=\int_{\GG\times E}1_A(\gg)\gg(\d z)\pi_\mu(\d\gg)\\
&=\int_A\gg(E)\pi_\mu(\d\gg)=0.
\end{aligned}
\end{eqnarray*}
\end{proof}

\beg{prp}\label{P3.3}    $(\Eg,\D(\Eg))$ is a conservative symmetric
Dirichlet form on $L^2(\pi_\mu)$ with $\D(\Eg)$ including the family
of cylindrical functions
\beg{equation*}\beg{split}
\F_\mu^C :=\Big\{\gg\mapsto & f\big(\gg(h_1),\cdots, \gg(h_m)\big): \ m\ge 1, f\in C_b^1(\R^m), \\
&\qquad\qquad\qquad h_i\in L^1(\mu)\cap L^\infty(\mu), \|\psi_\mu 1_{  h_i\ne 0}\|_\infty<\infty\Big\},
\end{split}\end{equation*}
where
$\|\cdot\|_\infty$ is the $L^\infty(\mu)$-norm. \end{prp}

 \beg{proof} According to Lemma \ref{L3.2}, for $F,G\in \D(\Eg)$,
 $\Eg(F,G)$ is finite and does not depend on $\pi_\mu$-versions of $F$ and $G$. Thus, $(\Eg,\D(\Eg))$ is a well defined positive bilinear form on $L^2(\pi_\mu)$.   Since  $\F_\mu^C$ is dense in $L^2(\pi_\mu)$ and the normal contractivity property is trivial by the definition of $\Eg$,    it remains to show   $\D(\Eg)\supset \F_\mu^C$ and the closed property of the form.   We prove
 these two points separately.

\medskip

\medskip

(a) Let $F\in \F_\mu^C$ with
$$F(\gg)=f\big(\gg(h_1),\cdots, \gg(h_m)\big),\quad\gg\in\Gamma,$$ which is well defined in $L^2(\pi_\mu)$
since $\gg(K)<\infty$ for $\pi_\mu$-a.e. $\gg\in\GG$ and any compact
subset $K$ of $E$. We intend to show that $\Eg(F,F)<\infty$. Since
$f\in C_0^1(\R^m)$, $h_i\in L^1(\mu)\cap L^\infty(\mu)$, and there exists $n\geq1$ such
that
$$\mu\big(h_i\neq0,\psi_\mu>n\big)=0,\quad i=1,\cdots,m,$$
we obtain
\beg{equation*}\beg{split}
\Eg(F,F)&=\int_{\GG\times(\bigcup_{i=1}^m \{h_i\neq0\})}\Big[f\big(\gg(h_1)+h_1(z),\cdots, \gg(h_m)+h_m(z)\big)\\
&\qquad\qquad\qquad\qquad\quad-f\big(\gg(h_1),\cdots, \gg(h_m)\big)\Big]^2
\psi(\gg,z)\pi_\mu(\d\gg)\mu(\d z)\\
&\le \|\nn f\|_\infty^2 \int_{\GG\times\{\psi_\mu\leq n\}}\sum_{i=1}^mh_i(z)^2\psi(\gg,z)\pi_\mu(\d\gg)\mu(\d z)\\
&=\|\nn f\|_\infty^2 \sum_{i=1}^m\int_{\{\psi_\mu\leq n\}}h_i(z)^2\psi_\mu(z)\mu(\d z)\\
&\le n\|\nn f\|_\infty^2\sum_{i=1}^m \mu(h_i^2)\leq n\|\nn f\|_\infty^2\sum_{i=1}^m\|h_i\|_\infty\mu(|h_i|)<\infty.
\end{split}\end{equation*}

\medskip

(b) Let $\{F_n\}_{n\geq1}$ be an $\Eg_1$-Cauchy sequence. We shall find
$F\in\D(\Eg)$ such that $\Eg_1(F_n-F,F_n-F):=\Eg(F_n-F,F_n-F)+\pi_\mu(|F_n-F|^2)\rightarrow0$ as
$n\rightarrow\infty$. Since $\{F_n\}_{n\geq1}$ is a Cauchy sequence in
$L^2(\pi_\mu)$ (which is complete), there exists $F\in L^2(\pi_\mu)$
such that $F_n\rightarrow F$ in $L^2(\pi_\mu)$. Now we can choose a
subsequence $\{F_{n_k}\}_{k\geq1}$ such that $F_{n_k}\rightarrow F$
$\pi_\mu$-a.e. By  Lemma \ref{L3.2} we have $F_{n_k}(\gg+\dd_z)\to
F(\gg+\dd_z)$ for $(\pi_\mu\times\mu)$-a.e. $(\gg,z)\in\GG\times E$. Therefore, it follows from the Fatou
lemma that
\begin{eqnarray*}
\begin{aligned}
&\Eg(F_n-F,F_n-F)\\
&=\int_{\Gamma\times
E}\liminf_{n_k\rightarrow\infty}\big[(F_n-F_{n_k})(\gg+\delta_z)-(F_n-F_{n_k})(\gg)\big]^2
\psi(\gg,z)\pi_\mu(\d\gg)\mu(\d z)\\
&\leq\liminf_{n_k\rightarrow\infty}\Eg(F_n-F_{n_k},F_n-F_{n_k}).
\end{aligned}
\end{eqnarray*}
Since $\{F_n\}_{n\geq1}$ is an $\Eg$-Cauchy sequence and $F_n\to F$ in $L^2(\pi_\mu)$, this implies that
$$\lim_{n\rightarrow\infty}\Eg_1(F_n-F,F_n-F)=0.$$ Combining this
with the fact that
$$\Eg(F,F)\leq2\Eg(F_n-F,F_n-F)+2\Eg(F_n,F_n),\quad n\geq1,$$
we conclude that $F\in\D(\Eg)$ and $F_n\rightarrow F$ in $\D(\Eg)$ as $n\rightarrow\infty$.
\end{proof}

The next result provides a criterion for the regularity of the
Dirichlet form, which ensures the existence of the associated Markov
process according to the Dirichlet form theory (see \cite{f, MR}). To this end, we first reduce $\GG$ to a locally compact subspace
$\GG_\mu$. Since $\GG$ is a Polish space such that the set $\{\pi_\mu\}$ of single probability measure is tight, we can
choose an increasing sequence $\{K_n\}_{n\ge 1}$ consisting of compact subsets of $\GG$ such that $\pi_\mu(K_n^c)\leq1/n$
for any $n\geq1$.
Then $\pi_\mu$ has full measure on $\GG_\mu:=\cup_{n=1}^\infty K_n$, which is a locally compact separable metric space.

\beg{prp} \label{P3.4}    If $\psi\in L^1(\pi_\mu\times\mu)$, then  $(\Eg,\D(\Eg))$ is a regular
Dirichlet form on $L^2(\GG_\mu; \pi_\mu)$. \end{prp}

\beg{proof} Since $\psi\in L^1(\pi_\mu\times\mu)$, we have
$\BB_b(\GG_\mu)\subset\D(\Eg)$, where $\BB_b(\GG_\mu)$ is the set of all
bounded measurable functions on $\GG_\mu$. In particular,
$C_0(\GG_\mu)\subset\D(\Eg)$. Thus, it suffices to prove that
$C_0(\GG_\mu)$ is dense in $\D(\Eg)$ w.r.t. the $\Eg_1$-norm, i.e. for
any $F\in\D(\Eg)$, one may find a sequence $\{F_n\}_{n\geq1}\subset
C_0(\GG_\mu)$ such that $\Eg_1(F_n-F,F_n-F)\rightarrow0$ as
$n\rightarrow\infty$.

\medskip

Since $\BB_b(\GG_\mu)\cap\D(\Eg)$ is dense in $\D(\Eg)$ (see e.g. \cite[Proposition
I.4.17]{MR}), we may assume that $F\in \BB_b(\GG_\mu)$. Moreover, since
$C_0(\GG_\mu)$ is dense in $L^2(\GG_\mu;\pi_\mu)$, we may find a sequence
$\{F_n\}_{n\geq1}\subset C_0(\GG_\mu)$ such that
$\sup_{n\in\mathbb{N}}\|F_n\|_\infty\le \|F\|_\infty$ and
$\pi_\mu(|F_n-F|^2)\to 0$ as $n\to\infty.$  Without loss of
generality, we assume furthermore that $F_n\rightarrow F$
$\pi_\mu$-a.e. By Lemma \ref{L3.2}, $F_n(\gg+\dd_z)\to F(\gg+\dd_z)$ and
$(F_n-F)^2(\gg+\delta_z)\leq(\|F_n\|_\infty+\|F\|_\infty)^2\leq4\|F\|_\infty^2$ for $(\pi_\mu\times\mu)$-a.e. $(\gg,z)\in\GG\times E$.

\medskip

Note that (we do not have to distinguish integrals on $\GG_\mu$ and $\GG$ since $\pi_\mu(\GG_\mu^c)=0$)
\begin{eqnarray*}
\begin{aligned}
&\Eg(F_n-F,F_n-F)\\
&\leq2\int_{\GG\times E}(F_n-F)^2(\gg+\dd_z)\psi(\gg,z)\pi_\mu(\d\gg)\mu(\d z)\\
&\quad+2\int_{\GG\times E}(F_n-F)^2(\gg)\psi(\gg,z)\pi_\mu(\d
\gg)\mu(\d z).
\end{aligned}
\end{eqnarray*}
Since $\psi\in L^1(\pi_\mu\times\mu)$, by the dominated convergence theorem we obtain
$$\lim_{n\to\infty} \Eg(F_n-F, F_n-F)=0.$$
Combining this with $\pi_\mu(|F_n-F|^2)\to 0$, we conclude that
$$\lim_{n\rightarrow\infty}\Eg_1(F_n-F,F_n-F)=0,$$
which completes the proof.
\end{proof}

Finally,  we consider  the generator $(\scr L^\psi ,\D(\scr
L^\psi))$ of the Dirichlet form $(\Eg,\D(\Eg))$. For a measurable function $F$ on
$\GG$, let
\begin{eqnarray*}
\begin{aligned}
  &   \scr L^\psi_bF(\gg)=\int_{E}\left(F(\gg+\delta_z)-F(\gg)\right)\psi(\gg,z)\mu(\d z),\\
 & \scr L^\psi_d F(\gg)=
\int_E1_{\{\gg\geq\delta_z\}}\big(F(\gg-\delta_z)-F(\gg)\big)\psi(\gg-\delta_z,z)\gg(\d z),\quad \gg\in\GG
\end{aligned}
\end{eqnarray*}
provided the integrals above exist.
\beg{prp} Suppose $F\in \D(\Eg)$ such that $\scr L^\psi_bF, \scr L^\psi_d F\in L^2(\pi_\mu)$. Then $F\in\D(\scr L^\psi)$
and $\scr L^\psi  F=\scr L^\psi_d F+\scr L^\psi_bF$. In particular, if $\mu$ is locally finite and
\begin{equation}\label{loc}
\int_\GG\psi(\gg,\cdot)^2\pi_\mu(\d\gg)\in L^1_{\text{loc}}(\mu),
\end{equation}
then
$$
\D(\scr L^\psi)\supset\Big\{\gg\mapsto  f\big(\gg(h_1),\cdots, \gg(h_m)\big): \ m\ge 1, f\in C_b^1(\R^m),
h_i\in C_0(E)\Big\}.
$$
\end{prp}
\beg{proof} (1) For any $F\in\D(\Eg)$ such that $\scr L_b^\psi F, \scr L_d^\psi F\in L^2(\pi_\mu)$, by the Mecke identity
(\ref{M}) for
$$H(\gg,z)= F(\gg)1_{\{\gg\geq\delta_z\}}\big(F(\gg-\dd_z)-F(\gg)\big)\psi(\gg-\dd_z,z),$$ we obtain
\begin{eqnarray*}
\begin{aligned}
&-\Eg(F,F)\\
&=\int_{\Gamma\times E}F(\gg)\big(F(\gg+\delta_z)-F(\gg)\big)\psi(\gg,z)\pi_\mu(\d \gg)\mu(\d z)\\
&\quad+\int_{\Gamma\times E}F(\gg+\delta_z)\big(F(\gg)-F(\gg+\delta_z)\big)\psi(\gg,z)\pi_\mu(\d \gg)\mu(\d z)\\
&=\int_{\Gamma\times E}F(\gg)\big(F(\gg+\delta_z)-F(\gg)\big)\psi(\gg,z)\pi_\mu(\d \gg)\mu(\d z)\\
&\quad+\int_{\Gamma\times E}F(\gg)1_{\{\gg\geq\delta_z\}}\big(F(\gg-\delta_z)-F(\gg)\big)\psi(\gg-\delta_z,\gg)\gg(\d z)\pi_\mu(\d \gg)\\
&=\int_{\Gamma}F(\gg)\big(\scr L^\psi_bF+\scr L^\psi_dF\big)(\gg)\pi_\mu(\d
\gg).
\end{aligned}
\end{eqnarray*}
Hence, the first assertion follows.

\bigskip

(2) Let
$$F(\gg)=f\big(\gg(h_1),\cdots,\gg(h_m)\big),\quad \gg\in\GG,$$
where $f\in C_b^1(\R^m)$, $h_i\in C_0(E)$ and $m\geq1$.
By  the Schwartz inequality we have
\begin{eqnarray*}
\begin{aligned}
&\int_{\GG\times E}\big(F(\gg+\delta_z)-F(\gg)\big)^2\psi(\gg,z)^2\pi_\mu(\d\gg)\mu(\d z)\\
&=\int_{\GG\times(\bigcup_{i=1}^m \supp h_i)}\Big[f\big(\gg(h_1)+h_1(z),\cdots, \gg(h_m)+h_m(z)\big)\\
&\qquad\qquad\qquad\qquad\quad-f\big(\gg(h_1),\cdots, \gg(h_m)\big)\Big]^2
\psi(\gg,z)^2\pi_\mu(\d\gg)\mu(\d z)
\end{aligned}
\end{eqnarray*}
\begin{eqnarray*}
\begin{aligned}
&\le \|\nn f\|_\infty^2 \sum_{i=1}^m\int_{\GG\times(\bigcup_{i=1}^m \supp h_i)}h_i(z)^2\psi(\gg,z)^2\pi_\mu(\d\gg)\mu(\d z)\\
&\leq\|\nn f\|_\infty^2 \left(\sum_{i=1}^m\|h_i\|_\infty^2\right)\int_{\GG\times(\bigcup_{i=1}^m \supp h_i)}\psi(\gg,z)^2\pi_\mu(\d\gg)\mu(\d z)\\
&<\infty,
\end{aligned}
\end{eqnarray*}
where the last step is due to (\ref{loc}). Then $\scr L_b^\psi F\in L^2(\pi_\mu)$ since
$$\|\scr L_b^\psi F\|_{L^2(\pi_\mu)}^2\leq\int_{\GG\times E}\big(F(\gg+\delta_z)-F(\gg)\big)^2\psi(\gg,z)^2\pi_\mu(\d\gg)\mu(\d z)<\infty.$$
On the other hand, using the Mecke identity (\ref{M}) for
$$H(\gg,z)=1_{\{\gg\geq\delta_z\}}\big(F(\gg-\delta_z)-F(\gg)\big)^2\psi(\gg-\delta_z,z)^2,$$
we arrive at
\begin{eqnarray*}
\begin{aligned}
\|\scr L_d^\psi F\|_{L^2(\pi_\mu)}^2&\leq\int_{\GG\times E}1_{\{\gg\geq\delta_z\}}\big(F(\gg-\delta_z)-F(\gg)\big)^2\psi(\gg-\delta_z,z)^2\gg(\d z)\pi_\mu(\d\gg)\\
&=\int_{\GG\times E}\big(F(\gg+\delta_z)-F(\gg)\big)^2\psi(\gg,z)^2\pi_\mu(\d\gg)\mu(\d z)<\infty.
\end{aligned}
\end{eqnarray*}
Consequently, $\scr L_d^\psi F\in L^2(\pi_\mu)$ and the proof is now completed according to the first assertion.
\end{proof}

\section{The path space of L\'evy processes}

Let $X=\{X_t:t\geq0\}$ be the L\'{e}vy process on $\R^d$ starting
from $0$ with a constant drift $b\in\R^d$ and the L\'{e}vy measure
$\nu$, which satisfies $\nu(\{0\})=0$ and
$$
\int_{\R^d}\left(|z|^2\wedge1\right)\nu(\d z)<\infty.
$$
So, $X_t$ is generated by
$$
\scr L f=\<b,\nn f\>+\int_{\R^d} \big\{f(z+\cdot)-f -\<\nn f,
z\>1_{\{|z|\le 1\}}\big\}\nu(\d z),
$$
which is well defined for $ f\in C_b^2(\R^d)$.

\medskip

\medskip

Let $\Lambda$ be the distribution of $X$, which is a probability
measure on the path space
$$
W:=\{w:[0,\infty)\rightarrow\R^d\,|\,w\,\,\mbox{is right continuous
having left limits}\}.
$$
It is well known that $W$ is a Polish space under the Skorokhod
metric
\begin{eqnarray*}
\begin{aligned}
\text{dist}(v,w):=&\inf\bigg\{\delta>0:\mbox{there
exist}\,\, n\geq1,\, 0=s_0<s_1<\cdots<s_n,\,\, \mbox{and}\,\,0=t_0<t_1\\
&\quad\quad\quad\quad<\cdots<t_n\quad\mbox{such that}\quad|t_i-s_i|\leq\delta\quad\mbox{and} \\
&\quad\qquad\quad\quad\quad\displaystyle\sup_{s\in
[s_{i-1},s_i),t\in[t_{i-1},t_i)}1\land
|v_s-w_t|\leq\delta\quad\mbox{hold for all}\quad 1\leq i\leq
n\bigg\}.
\end{aligned}
\end{eqnarray*}

\medskip

Let $\ttpsi\in L^1(\LL\times\nu\times\d t)$ be a nonnegative measurable function on $W\times(\R^d\setminus\{0\})\times[0,\infty)$ such that
 $$\ttpsi_{\nu\times\d t}(x,t):= \int_{W} \ttpsi(w,x,t)\LL(\d w)<\infty,\ \ \ \ (\nu\times\d t)\text{-a.e.}\  (x,t)\in (\R^d\setminus\{0\})\times[0,\infty).$$
Consider
\begin{eqnarray*}
\begin{aligned}
\Etg(F,G)&:=\int_{W\times\R^d\times[0,\infty)}\left(F(w+x1_{[t,\infty]})-F(w)\right)
\left(G(w+x1_{[t,\infty]})-G(w)\right)\\
&\quad\qquad\qquad\qquad\times\ttpsi(w,x,t)\,\Lambda(\d w)\nu(\d
x)\d t
\end{aligned}
\end{eqnarray*}
for
$$
F,\,G\in\D(\Etg):=\{F\in L^2(\Lambda):\Etg(F,F)<\infty\}.
$$

\medskip

\medskip

To apply the known Poincar\'e inequality on Poisson space,  we
follow the line of \cite{wu00} by constructing the L\'evy process
using  Poisson point processes. Let $E=( \R^d\setminus\{0\})\times
[0,\infty)$, which is a Polish space by taking the following
complete metric on $\R^d\setminus\{0\}$:
\begin{eqnarray*}
\begin{aligned}
\rho(x,y)&:=\sup\bigg\{|f(x)-f(y)|:|\nabla f(z)|\leq\frac1{|z|}\vee1,\\
&\qquad\qquad\qquad\qquad\qquad\qquad z\in\R^d\setminus\{0\},f\in
C^1\left(\R^d\setminus\{0\}\right)\bigg\}.
\end{aligned}
\end{eqnarray*}

Next, let $\mu= \nu\times \d t$, which is finite on bounded subsets of $E$ and
does not have atom. Let $\pi_\mu$ be the Poisson measure with
intensity $\mu$, which is a probability measure on the configuration
space
$$
\Gamma:=\left\{\sum_{i=1}^n\delta_{(x_i,t_i)}:x_i\in\R^d\setminus\{0\},t_i\in[0,\infty),1\leq
i\leq n, n\in\mathbb{Z}_+\cup\{\infty\}\right\}.
$$
Then on the probability space $(\Gamma,\scr F_\Gamma,\pi_\mu)$, the
L\'{e}vy process $X_t$ can be formulated as (see \cite{js})
$$
X_t(\gamma)=bt+\int_{\{|z|>1\}\times[0,t]}z\,\gamma(\d z,\d s)+
\int_{\{|z|\leq1\}\times[0,t]}z\,(\gamma-\mu)(\d z,\d s),\quad
t\geq0,
$$
where the second term in the right hand side above is the Stieltjes
integral, and the last term is the It\^{o} integral. Therefore,
\beq\label{relation} \Lambda=\pi_\mu\circ X^{-1}.
\end{equation}
Combining this with the Mecke identity (\ref{M}), we obtain
\beq\label{M'}  \beg{split}
&\int_{W}\sum_{\triangle w_t\neq0}h(w,\triangle w_t,t)\Lambda(\d w)\\
&\quad\qquad=\int_{W\times(\R^d\setminus\{0\})\times[0,\infty)}h(w+x1_{[t,T]},x,t)\Lambda(\d
w)\nu(\d x)\d t
\end{split}\end{equation}  for any non-negative measurable function $h$ on $W\times\R^d\times[0,\infty)$.
Due to (\ref{relation}) and (\ref{M'}), arguments used  in Section 3
also work  for $(\Etg,\D(\Etg)),\LL$ and $\ttpsi$ in place of
$(\Eg,\D(\Eg)), \pi_\mu$ and $\psi$ respectively. In particular,
letting $\tt\ll_{L}(\tt\psi), \tt\ll_{T}(\tt\psi)$ and
$\tt\ll_{E}(\tt\psi)$ be, respectively, the exponential convergence
rates in the $L^2$-norm, the $L^2$-tail norm and entropy for the semigroup
associated to $(\Etg,\D(\Etg))$, we obtain the following result.
\begin{thm}
We have
$$\ess_{\LL\times\mu}\inf \tt\psi\le \tt\ll_{L}(\tt\psi), \tt\ll_{T}(\tt\psi), \tt\ll_{E}(\tt\psi)
\le \ess_\mu\inf \tt\psi_\mu,$$
and the equalities hold provided $\tt\psi(w,x,t)$ does not depend
on $w$.
\end{thm}

\beg{thebibliography}{99}

\bibitem{f} M. Fukushima,  Y. Oshima and M. Takeda,
\emph{Dirichlet Forms and Symmetric Markov Processes,} Walter de
Gruyter, 1994.

\bibitem{G}
L. Gross, \emph{ Logarithmic Sobolev inequalities,} Amer. J. Math.
97(1975), 1061--1083.

\bibitem{js}
J. Jacod and A. N. Shiryaev, \emph{Limit Theorems for Stochastic
Processes}, Springer, (1987), Berlin.

\bibitem{MR} Z.-M. Ma  and M. R\"ockner, \emph{Introduction to
the Theory of (Non-Symmetric) Dirichlet Forms,} Springer-Verlag,
Berlin, 1992.

\bibitem{me}
J. Mecke, \emph{Stationaire Zuf\"{a}llige Ma$\beta$e auf
lokalkompakten abelschen Gruppen}, Z. Wahrsch. verw. Geb. 9(1967),
36-58.

\bibitem{ro}
M. R\"{o}ckner, \emph{Stochastic analysis on configuration spaces:
basic ideas and recent results}, in ``New Directions in Dirichlet
forms,''  157-231, AMS/IP Stud. Adv. Math.  8, Amer. Math. Soc. Providence, RI, 1998.

\bibitem{RW} M. R\"ockner and F.-Y. Wang, \emph{Weak Poincar\'e inequalities
and $L^2$-convergence rates of Markov semigroups,} J. Funct. Anal.
185(2001), 564--603.

\bibitem{Simon} B. Simon, \emph{The $P(\Phi)_2$-Euclidean (Quantum) Field Theory,}
Princeton Univ. Press, Princeton, NJ. 1974.

\bibitem{S}
D. Surgailis, \emph{On the multiple Poisson stochastic integrals and
associated Markov semigroups}, Prob. and Math. Stat.  3(1984):
217-239.

\bibitem{W02} F.-Y. Wang, \emph{Functional inequalities and spectrum
estimates: the infinite measure case,} J. Funct. Anal. 194(2002),
288--310.

\bibitem{w05}
F.-Y. Wang, \emph{Functional Inequalities, Markov Processes and
Spectral Theory}, Science Press,  2005.

\bibitem{wu00}
L. Wu, \emph{A new modified logarithmic Sobolev inequality for
Poisson point processes and several applications}, Probab. Theory
Relat. Fields, 118(2000): 427-438.


\bibitem{Wu} L. Wu, \emph{Uniform positive improvingness, tail norm condition and spectral
gap,} available online
http://www.math.kyoto-u.ac.jp/probability/sympo/sa01/wu.pd

\bibitem{MZ} S.-Y. Zhang and Y.-H. Mao, \emph{Exponential convergence rate in Boltzman-Shannon entropy},
Sci Sin, Ser A, 2000, 44(3): 280-285.

\end{thebibliography}
\end{document}